\documentclass{article}
\usepackage{amssymb}
\usepackage{amsmath}  
\usepackage{amsthm}  

\usepackage{amscd}
\usepackage[utf8]{inputenc}
\usepackage{amsfonts}
\usepackage{amstext}
\usepackage{color}
\newtheorem{theorem}{Theorem}
\newtheorem{lemma}{Lemma}

 \newtheorem{remark}{Remark}
 
\newcommand{\R}{{\mathbb R}}
\newcommand{\eps}{\varepsilon}
 \newcommand{\Z}{{\mathbb Z}} \newcommand{\N}{{\mathbb N}}
\newcommand{\T}{{\mathbb T}}

\newcommand{\supp}{{\rm supp\,}}

\newcommand\ack{\section*{Acknowledgement}}

\begin{document}
\author{Aleksei~Kulikov \quad Ilya~Zlotnikov }

\title{Contractive projections in Paley-Wiener spaces}
\maketitle
\begin{abstract}
Let $S_1$ and $S_2$ be disjoint finite unions of parallelepipeds. We describe necessary and sufficient conditions on the sets $S_1,S_2$ and exponents $p$ such that the canonical projection $P$ from $PW_{S_1\cup S_2}^p$ to $PW_{S_1}^p$ is a contraction.
\end{abstract}
\date{}\maketitle

\section{Introduction}

In the present short paper, we study the particular case of the following

\medskip

\noindent{\bf Main problem.}

Let $X,Y$ be spaces of functions. Assume that $Y$ is a subspace of $X$ and $ P:X \to Y$ is a projection. What assumptions should be imposed 
on $X, Y$, and $P$ to ensure that $P$ is a contraction?

\medskip

The case when $X$ is an $L^p$ space and $Y$ is an arbitrarily closed subspace of $L^p$ was completely solved by T. Andô \cite{A}. He showed that if $p\neq 2$ and $P$ leaves constants intact then $P$ is a conditional expectation with respect to some $\sigma$-algebra. He also obtained a complete characterization even without this assumption, see \cite{A}, Theorem 2.

One prominent example is the case $X = L^p(\T^d)$,\, $ Y = \{f\in L^p(\T^d)\mid \hat{f}(\lambda) = 0, \lambda \notin \Lambda\}$ and the projection $P$ being an idempotent Fourier multiplier. In this case the result of Andô has the following simple formulation.
\begin{theorem}
The projection $P$ is a contraction if and only if either $p = 2$ or $\Lambda\subset \Z^d$ is a coset.
\end{theorem}

Recently, O.F.~Brevig, J.~Ortega-Cerdà, and K.~Seip \cite{BOS} studied the contractivity of the similar idempotent Fourier multipliers in the case when $X$ is a Hardy space, that is $$X = H^p(\T^d) = \{f\in L^p(\T^d)\mid \hat{f}(n_1, n_2,\ldots, n_d) = 0 \text{  \rm if }\, n_k < 0 \,\text{\rm for some}\, k\}.$$
They showed that if $p\notin 2\N$ then the only contractions are the same as in the result of Andô, while for $p = 2k, k\in \N$ there exist non-trivial examples if $d\ge 3$. For the complete statement of their results, see \cite{BOS}, Theorem 1.2.

Since the Hardy spaces are the subsets of the Lebesgue spaces on the torus with some restrictions on the Fourier transform, it is natural to consider the analogue of the problem solved in \cite{BOS} in the setting of Euclidean spaces.
Specifically, in this note we consider the operators acting on the Paley-Wiener spaces. For the Fourier transform of the function $f$ from $L^p, 1\le p \le \infty,$ we fix the notation
$$\mathcal{F}(f)(t) = \int\limits_{\R^d} f(x)e^{-2 \pi i \langle x, t\rangle} dx, \quad t \in \R^d,$$
where $\langle \cdot, \cdot\rangle$ stands for the scalar product and the integral is understood in the sense of distributions.

We recall that for a compact set $S \subset \R^d$ and $1 \le p \le \infty$ the Paley-Wiener space $PW^p_S$ is defined by
$$
PW^p_S = \left\{ f \in L^p(\R^d) : {\rm Sp} (f):= \supp \mathcal{F}(f) \subset S \right\}.
$$
Sometimes, the spaces $PW^{\infty}_S$ are called  Bernstein spaces in the literature.
We refer the reader to books \cite{levin} and \cite{ou} for more details regarding properties of Paley-Wiener spaces.  

Assume that $S_1$ and $S_2$ are disjoint compact sets.
In what follows we deal with the {\it canonical projection} $P$ acting from $PW_{S_1\cup S_2}^p$ to $PW^p_{S_1}$ defined by
$$
P(f)(x):= \mathcal{F}^{-1}[\mathcal{F}(f)\cdot\chi_{S_1})](x).
$$

For a set $S\subset \R^d$ and $k\in \mathbb{N}_0$ we define $kS$ inductively as $0S = \{ 0\}$, $(k+1)S = kS + S$, where $+$ denotes the Minkowski sum.   
\begin{theorem}\label{main_th}
Let $S_1$ and $S_2$ be  finite unions of parallelepipeds in $\R^d$ such that  ${\rm mes}(S_1\cap S_2) = 0$ and let $P$ be a canonical projection from $PW_{S_1\cup S_2}^p$ to $PW^p_{S_1}$. We have
\begin{enumerate}

\item[\rm{1.}] If $p \in 2\N$ then $P$ is a contraction if and only if 
\begin{equation}\label{cond1}
{\rm mes}\left(\left(\frac{p}{2}S_1 + \left(\frac{p}{2}-1\right)(-S_1) \right) \cap S_2 \right) = 0.
\end{equation}
\item[\rm{2.}] If $p \notin 2\N$ then $P$ is a contraction if and only if ${\rm mes}(S_1) = 0$ or  \hbox{${\rm mes}(S_2) = 0$}.
\end{enumerate}
\end{theorem}
\begin{remark}
Note that from this theorem it follows that if $P$ is a contraction for $p = 2(n+k),\, n, k\in \N$ then it is a contraction for $p = 2n$ as well. On the other hand, for each $n\in \N$ there are sets $S_1, S_2$ such that $P$ is a contraction if and only if $p = 2m, m\in\N, m\le n$. For example, we could take $S_1 = [0, 1], S_2 = [n, n+1]$. Note also that we can not have projections for all $p\in 2\N$ {\rm(}unless ${\rm mes}(S_1) = 0$ or ${\rm mes}(S_2) = 0${\rm)} since for large enough $n$ we have
$$
{\rm mes}\left(\left(nS_1 + \left(n-1\right)(-S_1) \right) \cap S_2 \right) > 0.
$$

\end{remark}

\begin{remark}
We also note that, in contrast with the Hardy space setting, our argument does not depend on the dimension under the assumption that $S_1$ and $S_2$ are finite unions of parallelepipeds {\rm(}in  fact, the proof is still valid under even more general assumptions, see Remark~{\rm\ref{cantor_rem}}{\rm)}.
\end{remark}

\section{Proofs} For simplicity, we consider the case $d=1$, i.e. the sets involved are disjoint unions of intervals. The proof of the general case is done exactly in the same way and therefore is omitted.

We follow the argument from \cite{BOS} (see Lemma~3.1 there) and invoke the criterion due to H.S.~Shapiro from \cite{Sh}.
\begin{lemma}\label{criterion}
Assume that $S_1$ and $S_2$ are disjoint compact sets.
Let $f \in PW^p_{S_1}$.
The following statements are equivalent:

\begin{enumerate}
    \item[\rm{(}a\rm{)}]  The inequality
    $$\|f\|_{L^p} \le \|f + g\|_{L^p}$$
    is true for every $g \in PW^p_{S_2}$.
    \item[\rm{(}b\rm{)}] The equality
    \begin{equation}\label{p_eq}
        \int\limits_{\R} |f(x)|^{p-2} f(x) \overline{g(x)} \,d x = 0
    \end{equation}
    holds for every $g \in PW^p_{S_2}$. 
\end{enumerate} 
\end{lemma}
For $1 < p < \infty$ this lemma immediately follows from \cite{Sh}, Theorem~4.2.1. In the case $p = 1$, the statement follows from \cite{Sh}, Theorem~4.2.2, since the set $\{x: f(x) = 0\}$ has measure zero for any non-trivial entire function $f$.
Informally, this lemma corresponds to taking the derivative of $||f+\eps g||_{L^p}^p$ at $\eps = 0$.
\subsection{\bf Proof of Theorem~\ref{main_th}, Part 1.}

Now, we consider $p \in 2\N$, i.e. $p =2k, k \in \N$. 
First, we note that the sufficiency of  condition~\eqref{cond1} follows from Lemma~\ref{criterion}. Indeed, take any $q \in PW^p_S$ and denote by $f = P(q)$ and $g = q -f$. Note that $f \in PW^p_{S_1}$ and $g \in PW^p_{S_2}$. Since Schwartz functions are dense in $PW^p_S$ if $S$ is a finite union of intervals, we will assume that $f, g$ are Schwartz functions so that all Fourier transforms are genuine functions.
By Titchmarsh's theorem, we have
$$
{\rm Sp}(|f(x)|^{p-2} f(x)) = {\rm Sp}( f^{k}(x) (\overline{f(x)})^{k-1}) \subset kS_1 + (k-1)(-S_1).$$
Using condition~\eqref{cond1} and Plancherel's theorem, we get
$$ \int\limits_{\R} |f(x)|^{p-2} f(x) \overline{g(x)} \,d x = \int\limits_{\R} \mathcal{F}(|f|^{p-2} \bar{f})(t) \overline{\mathcal{F}(g)(t)}\, dt = 0.
$$
Applying Lemma~\ref{criterion}, we finish the proof of the sufficiency:
$$
\|q\|_{L^p} = \|f + g\|_{L^p} \ge \|f\|_{L^p}.
$$

Second, we assume that for the sets $S_1$ and $S_2$  equation \eqref{cond1} does not hold and prove that the projection $P:PW^p_{S_1 \cup S_2} \to PW^p_{S_1},$ is not a contraction. We use again Lemma~\ref{criterion} and reduce the problem to the construction of  functions $f$ and $g$ that violate condition \eqref{p_eq}. Recall that $p = 2k, k \in \N,$ and set
$$
f(x) = \mathcal{F}(\chi_{S_1})(x),\quad g(x) = \mathcal{F}(\chi_{S_2})(x), 
$$
where $\chi_{S}$ stands for the indicator function of the set $S$. Note that $f, g\in L^p(\R)$ since $p \ge 2$.
Denote by
$$\Phi(t) = \mathcal{F}^{-1} (|f|^{p-2} f)(t).$$
Since $S_1$ is a finite union of intervals, ${\rm Sp}(\Phi) = kS_1 + (k-1)S_1$ and, moreover, $\Phi(t) > 0$ on the interior of $kS_1+(k-1)S_1$. We have 
$$
\int\limits_{\R}|f(x)|^{p-2} f(x) \overline{g(x)} \,d x = \int\limits_{S_2} \Phi(t) dt > 0.
$$
Thus, we arrive at a contradiction, since by Lemma~\ref{criterion} we have 
$$
\|f\|_{L^p} > \|f + g\|_{L^p}.
$$

\subsection{\bf Proof of Theorem~\ref{main_th}, Part 2.}

First, we prove that the canonical projection \hbox{$P:PW^p_{S_1\cup S_2} \to PW^p_{S_1}$} is not a contraction if $p$ is not an even integer and $1 \le p < \infty$.

We have to find functions $f$ from $PW^p_{S_1}$ and $g$ from $PW^p_{S_2}$ such that \eqref{p_eq} does not hold.
Clearly, it suffices to prove the following 
\begin{lemma}
Let $1\le p<\infty$ and $p \notin 2 \N$. Assume that $I$ and $J$ are non-empty disjoint intervals.
There is a function $f$ from the Paley-Wiener space $PW^p_I$ such that ${\rm mes}({\rm Sp}(|f|^{p-2}f) \cap J) > 0$.
\end{lemma}
\begin{proof}
Without loss of generality we can assume that $I = [-1, 1]$. First, we consider the case $p > 1$. Set 
\begin{equation}
    h(x) = (x^2 - 1)\frac{\cos(2\pi x) - \cosh(2\pi)}{\prod\limits_{k=1}^N \left((x+k)^2 + 1\right)},
\end{equation}
where $N > \frac{1}{2(p-1)} + 2$ so that $|h(x)|^{p-1}\in L^1(\R)$. Therefore, the function $g(x) = |h(x)|^{p-2}h(x)$ belongs to $L^1(\R)$ as well. Note that by the Paley-Wiener Theorem $h\in PW^p_I$. On the other hand, since $h$ changes sign at the points $\pm 1$, the function $g$ is not analytic. Therefore, its Fourier transform $G = \mathcal{F}(g)$ is a continuous function which is not compactly supported. Also, note that since $h$ is even, $G$ is even as well.

Let $J = [a, b]$ and assume without loss of generality that $a > 1$. Since $G$ is not compactly supported, there exists $x_0\in \R$ such that $x_0 > a$ and $x_0\in \supp\, G$ (here we used that $G$ is even so its support must extend both to $+\infty$ and $-\infty$). Since $G$ is continuous and $G(x_0)\neq 0$, we have $G(x)\neq 0, x_0 - \eps < x < x_0+\eps$ for small enough $\eps$. 

Consider $f(x) = h(\frac{a}{x_0}x)$. Since $\frac{a}{x_0} < 1$, $f$ also belongs to $PW^p_I$. On the other hand, we have $$ \mathcal{F}(|f|^{p-2}f)(\xi)= \frac{x_0}{a} G\left(\frac{x_0}{a}\xi\right).$$ Thus, $\mathcal{F}(|f|^{p-2}f)(a) = \frac{x_0}{a} G(x_0)\neq 0$ and ${\rm mes}(\supp\, \mathcal{F}(|f|^{p-2}f) \cap J) > 0$ as required.


\begin{remark}
Although it seems plausible that ${\rm Sp}\, G = \R$ we only managed to prove that ${\rm Sp}\, G$ is unbounded.
\end{remark}
Next, we deal with $p=1$. Consider the same function $h(x)$ with $N=2$. Note that $Q(x) = |f(x)|^{p-2}f(x) = 2\chi_{[-1,1]}(x) - 1$, whence 
$$
\mathcal{F}(Q)(x) = \frac{2 \sin(2\pi x)}{\pi x} - \delta_0,
$$
where $\delta_0$ is a Dirac delta measure at $0$. Clearly, the support of the distribution $\mathcal{F}(Q)$ is $\R$. This finishes the proof of the lemma.
\end{proof}
 
It remains to consider $p = \infty$. In this case we will show that the projection is not contractive directly. Again, we assume that $S_1 \supset [-1,1]$. Set
$$
f(x) = \frac{\sin(2\pi x)}{2 \pi x}, \quad f \in PW^{\infty}_{S_1}.
$$
Put $g(x) = \mathcal{F}(\chi_{S_2})(x)$. Clearly, $g(0) > 0$ and $g\in PW^\infty_{S_2}$. Consider $f_\eps(x) = f(x)-\eps g(x)$. We are going to show that for sufficiently small positive $\eps$ we have $\|f_{\eps}\|_{L^{\infty}} < 1$. This will contradict contractivity and finish the proof of the theorem.

Note that $f(x), g(x) \to 0$ as $|x|\to \infty$. Since $f(0) = 1$ and $|f(x)| < 1$ when $x\neq 0$, for every $\delta > 0$ there exists $\eps > 0$ such that $|f(x)-\eps g(x)| < 1$ if $|x| > \delta$. Thus, it remains to consider $|x| < \delta$.

We have $1-x^2 \le f(x) \le 1$ if $|x| < \frac{1}{100}$. On the other hand, $g(x) = c_0 + O(x)$ for some $c_0 > 0$ if $x$ is close enough to $0$. Thus,

$$0 < {\rm Re} (f(x) - \eps g(x)) \le 1-\eps c_0 + C \eps |x|$$
and
$$ |{\rm Im} (f(x)-\eps g(x))| \le C\eps |x|$$
for some constant $C$ and small enough $x$. Therefore,
$$|f(x)-\eps g(x)| \le 1-\eps c_0 + 2C\eps |x|.$$
Choosing $\delta$ so that $\delta < \frac{1}{100}$ and $2C\delta < c_0$, we get that 
$$|f(x) - \eps g(x)| < 1$$
for $|x| \le \delta$. Thus, $||f-\eps g||_{L^\infty(\R)} < 1$ and the projection is not contractive.
\begin{remark}\label{cantor_rem}
In fact, our results hold more generally when the sets $S_1, S_2$ are closures of open sets with boundaries of measure $0$ and results for $p\notin 2\N$ require only that $S_1, S_2$ contain some open balls. On the other hand, it would be interesting to study the case of arbitrary closed sets $S_1, S_2$ of positive measure with empty interior.
\end{remark}

\ack{The authors are grateful to Ole Fredrik Brevig, Joaquim Ortega-Cerdà, and Kristian Seip  for valuable discussions and constructive remarks.}

\noindent Aleksei Kulikov\\
Norwegian University of Science and Technology, Department of Mathematical Sciences\\
NO-7491 Trondheim, Norway\\
lyosha.kulikov@mail.ru

\bigskip

\noindent Ilya Zlotnikov\\
University of Stavanger, Department of Mathematics and Physics,\\
4036 Stavanger, Norway,\\
ilia.k.zlotnikov@uis.no

\end{document}